\documentclass[11pt]{article}
\usepackage[numbers,sort&compress]{natbib}
\usepackage[colorlinks,
            linkcolor=blue,
            anchorcolor=green,
            citecolor=magenta
            ]{hyperref}
\usepackage{mathrsfs,amssymb,amsfonts}
\usepackage{graphicx}
\usepackage{ifpdf}
\usepackage{amsmath,amssymb,latexsym,color}
\usepackage[mathscr]{eucal}
\usepackage[numbers]{natbib}
\newtheorem{thm}{Theorem}[section]

\newtheorem{prop}[thm]{Proposition}

\newtheorem{cor}[thm]{Corollary}
\newtheorem{example}{Example}
\newtheorem{remark}{Remark}[section]

\newcommand{\proof}{{\it Proof.\quad}}
\newcommand{\diam}{{\rm diam}}
\newcommand{\spec}{{\rm spec}}
\newcommand{\qed}{\hfill\Box\medskip}
\usepackage{CJK}
\begin{document}
\begin{CJK*}{GBK}{song}
\renewcommand{\abovewithdelims}[2]{
\genfrac{[}{]}{0pt}{}{#1}{#2}}

\title{\bf On the spectra of strong power graphs of finite groups}

\author{ Xuanlong Ma\footnote{{\em E-mail address:} xuanlma@mail.bnu.edu.cn.}
 \\
{\footnotesize   \em  Sch. Math. Sci. {\rm \&} Lab. Math. Com. Sys., Beijing Normal University, Beijing, 100875,  China
}
}
 \date{}
 \date{}
 \maketitle

\begin{abstract}
We give
the characteristic polynomial of the distance or adjacency matrix of the strong power graph of a finite group, and
compute its distance and adjacency spectrum.


\end{abstract}

\section{Introduction}

Given a connected graph $\Gamma$, denote by $V(\Gamma)$ and $E(\Gamma)$ the vertex set and edge set, respectively. Let $V(\Gamma)=\{v_1,v_2,\ldots,v_n\}$. The {\em distance} between the vertices $v_i$ and $v_j$, denoted by $d_{\Gamma}(v_i,v_j)$, is the length of the shortest path between them. The {\em diameter} of $\Gamma$, denoted by \diam($\Gamma$), is the maximum distance between any pair of
vertices of $\Gamma$. The set of neighbours of a vertex $v_i$ in $\Gamma$ is denoted by
$N_{\Gamma}(v_i)$, that is, $N_{\Gamma}(v_i)=\{v_j\in V(\Gamma): \{v_i,v_j\}\in E(\Gamma)\}$.

The {\em distance matrix} $D(\Gamma)$ of $\Gamma$ is the $n\times n$ matrix, indexed by $V(\Gamma)$,
such that $D(\Gamma)_{v_i,v_j}=d_{\Gamma}(v_i,v_j)$. The characteristic polynomial
$\Theta(\Gamma,x)=|xI-D(\Gamma)|$ is the {\em distance characteristic polynomial} of $\Gamma$. 
Note that $D(\Gamma)$ is symmetric.
The distance characteristic polynomial has real roots $\mu_1\ge \mu_2\ge \dots \ge \mu_n$.
If $\mu_{i_1}\ge \mu_{i_2}\ge \dots \ge \mu_{i_t}$ are the distinct roots of  $\Theta(\Gamma,x)$, then the {\em $D$-spectrum} of $\Gamma$ can be written as
$$
\spec_D(\Gamma)=\left(
                 \begin{array}{cccc}
                   \mu_{i_1} & \mu_{i_2} & \dots & \mu_{i_t} \\
                   m_1 & m_2 & \dots & m_t \\
                 \end{array}
               \right),
$$
where $m_j$ is the algebraic multiplicity of $\mu_{i_j}$.
Clearly $\sum_{j=1}^{t} m_j=n$.
The {\em adjacency matrix} $A(\Gamma)$ of $\Gamma$ is a $n\times n$ matrix and indexed by $V(\Gamma)$. The $ij$-th entry of $A(\Gamma)$ is $1$ if the vertices $v_i$ and $v_j$ are adjacent, otherwise it is $0$. Denote by $\Phi(\Gamma,x)$ the characteristic polynomial of $A(\Gamma)$.
Similarly, we define the adjacency spectrum $\spec(\Gamma)$ of $\Gamma$.
The largest root of $\Theta(\Gamma,x)$(resp. $\Phi(\Gamma,x)$) is called the {\em distance spectral radius}(resp. {\em adjacency spectral radius}) of $\Gamma$.


Let $G$ denote a finite group of order $n$.
The {\it power graph} of $G$ was introduced by Chakrabarty et al. \cite{CGS09}, which has the vertex set $G$ and two distinct elements are adjacent if one is a power of the other. Motivated by this, Singh and Manilal \cite{SM10} defined the strong power graphs as
a generalization of the power graphs. The {\em strong power graph} $\mathcal{P}_s(G)$ of $G$ is a graph whose vertex set consists of the elements of $G$ and two distinct vertices $x$ and $y$ are adjacent if $x^{n_1}=y^{n_2}$ for some positive integers $n_1,n_2 < n$.
Very recently, Bhuniya and Berathe \cite{BB15} investigated the Laplacian of strong power graphs.

In this paper we study  
the characteristic polynomial of the distance or adjacency matrix of the strong power graph of a finite group, and
compute its distance and adjacency spectrum.

\section{The results}\label{main}
Throughout this section $G$ denotes a finite group, and $\mathbb{Z}_n$ stands for the cyclic group of order $n$. We always assume $\mathbb{Z}_n=\{0,1,\ldots,n-1\}$.
For strong power graphs we have the following proposition.
\begin{prop}\label{basic}
(1) If $G$ is not cyclic, then $\mathcal{P}_s(G)$ is complete.

(2) $\mathcal{P}_s(\mathbb{Z}_n)$ is not connected if and only if $n$ is a prime number.

(3) $N_{\mathcal{P}_s(\mathbb{Z}_n)}(0)=\{
k\in \mathbb{Z}_n: m\ne 0, (m,n)\ne 1\}$, and the subgraph of $\mathcal{P}_s(\mathbb{Z}_n)$ induced by $\mathbb{Z}_n\setminus \{0\}$ is complete. In particular, \diam$(\mathcal{P}_s(\mathbb{Z}_n))=2$ if $n$ is not a prime number.
\end{prop}

Now we find the characteristic polynomial of the distance matrix  associated with the strong power graph $\mathcal{P}_s(\mathbb{Z}_n)$ for any composite number $n$.
\begin{thm}\label{thm1}
For any composite number $n$,
$$\Theta(\mathcal{P}_s(\mathbb{Z}_n),x)=(x+1)^{n-3}\Big(x^3+(3-n)x^2+(3-2n-3\phi(n))x
-\phi(n)^2-\phi(n)(4-n)-n+1\Big),$$
where $\phi(n)$ is Euler's totient function.
\end{thm}
\proof Write $\phi(n)=t$ and $k=n-\phi(n)-1$.
By Proposition \ref{basic},
the distance matrix $D(\mathcal{P}_s(\mathbb{Z}_n))$ is the $n\times n$ matrix given below,
where the rows and columns are indexed in order by the vertices in $N_{\mathcal{P}_s(\mathbb{Z}_n)}(0)$ then all generator elements of $\mathbb{Z}_n$, and $0$ is in last position.
\begin{equation*}\label{}
\setlength{\arraycolsep}{8pt}
\renewcommand\arraystretch{0.9}
D(\mathcal{P}_s(\mathbb{Z}_n))=
\left(
  \begin{array}{cccccccc}
    0 & 1 & \dots & 1 & 1 & \dots & 1 & 1 \\
    1 & 0 & \dots & 1 & 1 & \ldots & 1 & 1 \\
    \vdots & \vdots & \ddots & \vdots & \vdots & \ddots & \vdots & \vdots \\
    1 & 1 & \ldots & 0 & 1 & \ldots & 1 & 1 \\
    1 & 1 & \ldots & 1 & 0 & \ldots & 1 & 2 \\
    \vdots & \vdots & \ddots & \vdots & \vdots & \ddots & \vdots & \vdots \\
    1 & 1 & \ldots & 1 & 1 & \ldots & 0 & 2 \\
    1 & 1 & \ldots & 1 & 2 & \ldots & 2 & 0 \\
  \end{array}
\right).
\end{equation*}
The characteristic polynomial of $D(\mathcal{P}_s(\mathbb{Z}_n))$ is
\begin{equation}\label{de1}
\setlength{\arraycolsep}{8pt}
\renewcommand\arraystretch{0.9}
\Theta(\mathcal{P}_s(\mathbb{Z}_n),x)=
\left|
  \begin{array}{cccccccc}
    x & -1 & \dots & -1 & -1 & \dots & -1 & -1 \\
    -1 & x & \dots & -1 & -1 & \ldots & -1 & -1 \\
    \vdots & \vdots & \ddots & \vdots & \vdots & \ddots & \vdots & \vdots \\
    -1 & -1 & \ldots & x & -1 & \ldots & -1 & -1 \\
    -1 & -1 & \ldots & -1 & x & \ldots & -1 & -2 \\
    \vdots & \vdots & \ddots & \vdots & \vdots & \ddots & \vdots & \vdots \\
    -1 & -1 & \ldots & -1 & -1 & \ldots & x & -2 \\
    -1 & -1 & \ldots & -1 & -2 & \ldots & -2 & x \\
  \end{array}
\right|.
\end{equation}
Subtract the first column from the columns $2,3,\ldots,n$ of (\ref{de1}) to obtain (\ref{de2}):
\begin{equation}\label{de2}
\setlength{\arraycolsep}{8pt}
\renewcommand\arraystretch{0.9}
(x+1)^{k-1}
\left|
  \begin{array}{cccccccc}
    x & -1 & \dots & -1 & -1-x & \dots & -1-x & -1-x \\
    -1 & 1 & \dots & 0 & 0 & \ldots & 0 & 0 \\
    \vdots & \vdots & \ddots & \vdots & \vdots & \ddots & \vdots & \vdots \\
    -1 & 0 & \ldots & 1 & 0 & \ldots & 0 & 0 \\
    -1 & 0 & \ldots & 0 & x+1 & \ldots & 0 & -1 \\
    \vdots & \vdots & \ddots & \vdots & \vdots & \ddots & \vdots & \vdots \\
    -1 & 0 & \ldots & 0 & 0 & \ldots & x+1 & -1 \\
    -1 & 0 & \ldots & 0 & -1 & \ldots & -1 & x+1 \\
  \end{array}
\right|.
\end{equation}
Adding the rows $2,3,\ldots,n-1$ to the first row of (\ref{de2}). Then adding columns $2,3,\ldots,k$ to the first column,
we arrive at the determinant (\ref{de3}):
\begin{equation}\label{de3}
\setlength{\arraycolsep}{8pt}
\renewcommand\arraystretch{0.9}
(x+1)^{k-1}
\left|
  \begin{array}{cccccccc}
    x-n+2 & 0 & \dots & 0 & 0 & \dots & 0 & -1-x-t \\
    0 & 1 & \dots & 0 & 0 & \ldots & 0 & 0 \\
    \vdots & \vdots & \ddots & \vdots & \vdots & \ddots & \vdots & \vdots \\
    0 & 0 & \ldots & 1 & 0 & \ldots & 0 & 0 \\
    -1 & 0 & \ldots & 0 & x+1 & \ldots & 0 & -1 \\
    \vdots & \vdots & \ddots & \vdots & \vdots & \ddots & \vdots & \vdots \\
    -1 & 0 & \ldots & 0 & 0 & \ldots & x+1 & -1 \\
    -1 & 0 & \ldots & 0 & -1 & \ldots & -1 & x+1 \\
  \end{array}
\right|.
\end{equation}
Subtract the first column from the last column of (\ref{de3}). Then subtract the row
$k+1$ from the rows $k+2,k+3,\ldots,n$ to obtain (\ref{de4}):
\begin{equation}\label{de4}
\setlength{\arraycolsep}{5pt}
\renewcommand\arraystretch{0.9}
(x+1)^{n-3}
\left|
  \begin{array}{ccccccccc}
    x-n+2 & 0 & \dots & 0 & 0 &0& \dots & 0 & -2x+k-2 \\
    0 & 1 & \dots & 0 & 0 & 0&\ldots & 0 & 0 \\
    \vdots & \vdots & \ddots & \vdots & \vdots &\vdots & \ddots & \vdots & \vdots \\
    0 & 0 & \ldots & 1 & 0 &0& \ldots & 0 & 0 \\
    -1 & 0 & \ldots & 0 & x+1 &0& \ldots & 0 & 0 \\
    0 & 0 & \ldots & 0 & -1 &1& \ldots & 0 & 0 \\
    \vdots & \vdots & \ddots & \vdots & \vdots &\vdots & \ddots & \vdots & \vdots  \\
    0 & 0 & \ldots & 0 & -1 & 0&\ldots & 1 & 0 \\
    0 & 0 & \ldots & 0 & -2-x &-1& \ldots & -1 & x+2 \\
  \end{array}
\right|.
\end{equation}
Adding the rows $k+2,k+3,\ldots,n-1$ to the last row of (\ref{de4}), we get (\ref{de5}):
\begin{equation}\label{de5}
\setlength{\arraycolsep}{5pt}
\renewcommand\arraystretch{0.9}
(x+1)^{n-3}
\left|
  \begin{array}{ccccccccc}
    x-n+2 & 0 & \dots & 0 & 0 &0& \dots & 0 & -2x+k-2 \\
    0 & 1 & \dots & 0 & 0 & 0&\ldots & 0 & 0 \\
    \vdots & \vdots & \ddots & \vdots & \vdots &\vdots & \ddots & \vdots & \vdots \\
    0 & 0 & \ldots & 1 & 0 &0& \ldots & 0 & 0 \\
    -1 & 0 & \ldots & 0 & x+1 &0& \ldots & 0 & 0 \\
    0 & 0 & \ldots & 0 & -1 &1& \ldots & 0 & 0 \\
    \vdots & \vdots & \ddots & \vdots & \vdots &\vdots & \ddots & \vdots & \vdots  \\
    0 & 0 & \ldots & 0 & -1 & 0&\ldots & 1 & 0 \\
    0 & 0 & \ldots & 0 & -x-t-1 &0& \ldots & 0 & x+2 \\
  \end{array}
\right|.
\end{equation}
Expand it along the first row to obtain (\ref{de6}):
\begin{equation}\label{de6}
(x+1)^{n-3}\Big((x-n+2)(x+1)(x+2)+(-1)^{1+n}(-2x+k-2)|A|\Big),
\end{equation}
where
\begin{equation}\label{de7}
\setlength{\arraycolsep}{7pt}
\renewcommand\arraystretch{0.9}
|A|=\left|
  \begin{array}{cccccccc}
    0 & 1 & \dots & 0 & 0 & 0&\ldots & 0  \\
    \vdots & \vdots & \ddots & \vdots & \vdots &\vdots & \ddots & \vdots  \\
    0 & 0 & \ldots & 1 & 0 &0& \ldots & 0  \\
    -1 & 0 & \ldots & 0 & x+1 &0& \ldots & 0  \\
    0 & 0 & \ldots & 0 & -1 &1& \ldots & 0  \\
    \vdots & \vdots & \ddots & \vdots & \vdots &\vdots & \ddots & \vdots   \\
    0 & 0 & \ldots & 0 & -1 & 0&\ldots & 1  \\
    0 & 0 & \ldots & 0 & -x-t-1 &0& \ldots & 0  \\
  \end{array}
\right|.
\end{equation}
In (\ref{de7}) the determinant has order $n-1$. Interchange the row $k-1$ and the last row of (\ref{de7}). Then expand it along the first column to obtain $|A|$ as follows:
\begin{equation*}\label{}
\setlength{\arraycolsep}{7pt}
\renewcommand\arraystretch{0.9}
|A|=(-1)^{n}(-x-t-1).
\end{equation*}
It follows that
$$
\Theta(\mathcal{P}_s(\mathbb{Z}_n),x)=(x+1)^{n-3}\Big(x^3+(3-n)x^2+(3-2n-3t)x-n+nt-t^2-4t+1\Big).
$$
This completes our proof. $\qed$

By Proposition \ref{basic}, one has that $\mathcal{P}_s(G)$ is connected if and only if
$G$ is not cyclic, or $G\cong \mathbb{Z}_n$ for some composite number $n$. Now we may  compute the $D$-spectrum of any connected strong power graph.

\begin{thm}\label{thm2}
If $G$ is not cyclic, then
$$
\spec_D(\mathcal{P}_s(G))=\left(
                            \begin{array}{cc}
                              n-1 & -1 \\
                              1 & n-1 \\
                            \end{array}
                          \right).
$$
If $G\cong \mathbb{Z}_n$ for some composite number $n$, then
$\spec_D(\mathcal{P}_s(\mathbb{Z}_n))$ is
$$
\setlength{\arraycolsep}{4pt}
\renewcommand\arraystretch{1.2}
\left(
                 \begin{array}{cccc}
                   -1 & \frac{n-3+2\cos\frac{\theta}{3}\sqrt{n^2+9\phi(n)}}{3} & \frac{n-3+2\cos\frac{\theta+2\pi}{3}\sqrt{n^2+9\phi(n)}}{3} & \frac{n-3+2\cos\frac{\theta-2\pi}{3}\sqrt{n^2+9\phi(n)}}{3} \\
                   n-3 & 1 & 1 & 1 \\
                 \end{array}
               \right),
$$
where $0<\theta <\frac{\pi}{2}$ and  $\theta=\arccos\frac{2n^3+27\phi(n)^2+27\phi(n)}{2\sqrt{(n^2+9\phi(n))^3}}$.
\end{thm}
\proof Note that if $G$ is not cyclic then $\mathcal{P}_s(G)$ is complete. Thus, it suffices to compute the distance spectrum of $\mathcal{P}_s(\mathbb{Z}_n)$ for  some composite number $n$.
Let $$f(x)=x^3+(3-n)x^2+(3-2n-3\phi(n))x
-\phi(n)^2-\phi(n)(4-n)-n+1.$$
Suppose that $f(-1)=0$. Then $\phi(n)(n-\phi(n)-1)=0$. It follows that $\phi(n)=n-1$.
Namely $n$ is a prime number, a contradiction. This implies that
$D$-spectrum of $\mathcal{P}_s(\mathbb{Z}_n)$ has $-1$ with multiplicity $n-3$ by Theorem \ref{thm1}.
Note that the canonical solutions of any quadratic and cubic equation. We conclude that
$f(x)$ has three pairwise distinct roots, as presented above. $\qed$

\begin{cor}\label{coro3}
For any composite number $n$, the distance spectral radius of $\mathcal{P}_s(\mathbb{Z}_n)$ is $$\frac{n-3+2\cos\frac{\theta}{3}\sqrt{n^2+9\phi(n)}}{3},$$
where $0<\theta <\frac{\pi}{2}$ and $\theta=\arccos\frac{2n^3+27\phi(n)^2+27\phi(n)}{2\sqrt{(n^2+9\phi(n))^3}}$.
\end{cor}

For any positive integer $n$, the adjacency matrix $A(\mathcal{P}_s(\mathbb{Z}_n))$ is given below,
where the rows and columns are indexed in order by the vertices in $N_{\mathcal{P}_s(\mathbb{Z}_n)}(0)$ then all generator elements of $\mathbb{Z}_n$, and $0$ is in last position. Note that $N_{\mathcal{P}_s(\mathbb{Z}_n)}(0)$ may be the empty set.
\begin{equation*}\label{}
\setlength{\arraycolsep}{8pt}
\renewcommand\arraystretch{0.9}
A(\mathcal{P}_s(\mathbb{Z}_n))=
\left(
  \begin{array}{cccccccc}
    0 & 1 & \dots & 1 & 1 & \dots & 1 & 1 \\
    1 & 0 & \dots & 1 & 1 & \ldots & 1 & 1 \\
    \vdots & \vdots & \ddots & \vdots & \vdots & \ddots & \vdots & \vdots \\
    1 & 1 & \ldots & 0 & 1 & \ldots & 1 & 1 \\
    1 & 1 & \ldots & 1 & 0 & \ldots & 1 & 0 \\
    \vdots & \vdots & \ddots & \vdots & \vdots & \ddots & \vdots & \vdots \\
    1 & 1 & \ldots & 1 & 1 & \ldots & 0 & 0 \\
    1 & 1 & \ldots & 1 & 0 & \ldots & 0 & 0 \\
  \end{array}
\right).
\end{equation*}

An argument similar to the one used in the computing of  $\Theta(\mathcal{P}_s(\mathbb{Z}_n),x)$,
we get $\Phi(\mathcal{P}_s(\mathbb{Z}_n),x)$.
\begin{thm}\label{thm3}
For any positive integer $n$,
$$\Phi(\mathcal{P}_s(\mathbb{Z}_n),x)=(x+1)^{n-3}\Big(x^3+(3-n)x^2+(3-2n+\phi(n))x
+(n-\phi(n)-1)(\phi(n)-1)\Big).$$
\end{thm}

\begin{cor}\label{}
For a prime number $p\ge 2$,
$$\Phi(\mathcal{P}_s(\mathbb{Z}_p),x)=x(x+1)^{p-2}(x+2-p).$$
\end{cor}

As an application of Theorem \ref{thm3}, we may obtain the  adjacency spectrum of the strong power graph of a finite group.

\begin{thm}\label{}
If $G$ is not cyclic, then
$$
\spec(\mathcal{P}_s(G))=\left(
                            \begin{array}{cc}
                              n-1 & -1 \\
                              1 & n-1 \\
                            \end{array}
                          \right).
$$
If $G\cong \mathbb{Z}_p$ for some prime number $p$, then
$$\spec(\mathcal{P}_s(\mathbb{Z}_n))=\left(
                                       \begin{array}{ccc}
                                         0 & -1 & p-2 \\
                                         - & p-2 & 1 \\
                                       \end{array}
                                     \right).
$$
If $G\cong \mathbb{Z}_n$ for some composite number $n$, then
$\spec(\mathcal{P}_s(\mathbb{Z}_n))$ is
$$
\setlength{\arraycolsep}{4pt}
\renewcommand\arraystretch{1.2}
\left(
                 \begin{array}{cccc}
                   -1 & \frac{n-3+2\cos\frac{\theta}{3}\sqrt{n^2-3\phi(n)}}{3} & \frac{n-3+2\cos\frac{\theta+2\pi}{3}\sqrt{n^2-3\phi(n)}}{3} & \frac{n-3+2\cos\frac{\theta-2\pi}{3}\sqrt{n^2-3\phi(n)}}{3} \\
                   n-3 & 1 & 1 & 1 \\
                 \end{array}
               \right),
$$
where $0<\theta <\frac{\pi}{2}$ and  $\theta=\arccos\frac{2n^3+27\phi(n)^2+27\phi(n)-36n\phi(n)}{2\sqrt{(n^2-3\phi(n))^3}}$.
\end{thm}

\begin{cor}\label{}
For any composite number $n$, the adjacency spectral radius of $\mathcal{P}_s(\mathbb{Z}_n)$ is $$\frac{n-3+2\cos\frac{\theta}{3}\sqrt{n^2-3\phi(n)}}{3},$$
where $0<\theta <\frac{\pi}{2}$ and $\theta=\arccos\frac{2n^3+27\phi(n)^2+27\phi(n)-36n\phi(n)}{2\sqrt{(n^2-3\phi(n))^3}}$.
\end{cor}

\section*{Acknowledgement}
This research is supported by National Natural Science Foundation of China (11271047,
11371204).

\end{CJK*}

\end{document}